\documentclass[11pt,a4paper]{article}
\usepackage{amsmath, amsfonts, amssymb, amsthm}
\usepackage{mathtools}
\usepackage{commath}               
\usepackage{paralist}
\usepackage{graphicx}  
\usepackage[T1]{fontenc}
\usepackage{float}

\usepackage{multirow}
\usepackage{algorithmic}
\usepackage{algorithm}
\usepackage{float}
\usepackage{wrapfig}
\usepackage{mdframed}
\usepackage{xcolor}
\usepackage{bm}
\usepackage[title]{appendix}
\usepackage{resizegather}
\usepackage{bbm}
\usepackage{subcaption} 
\usepackage{wrapfig}
\usepackage{epsfig}     
\usepackage{siunitx}
\usepackage{nicefrac}
\usepackage{xurl}
\usepackage{enumitem}
\usepackage{geometry}
\usepackage{tikz-cd} 
\usepackage{hyperref}  
\hypersetup{colorlinks=true, urlcolor=black, citecolor=red}

\usepackage{color}
\usepackage{todonotes}
\usepackage[normalem]{ulem}
\newcommand{\nir}[1]{{\color{black}{#1}}}
\newcommand{\vin}[1]{{\color{black}{#1}}}

\newcommand{\ns}[1]{{\color{black}{#1}}}

\newtheorem{theorem}{Theorem}

\newtheorem{lemma}[theorem]{Lemma}

\theoremstyle{definition}

\newtheorem{remark}{Remark}

\newcommand{\sumd}[1]{\sideset{}{'}\sum_{n=0}^{#1}}
\newcommand{\ppp}[1]{{ \text{P}_{#1} }}
\newcommand{\rrr}[2]{{ \text{R}_{ {#1} , {#2} }  }}
\newcommand{\cond}[1]{\operatorname{cond}\left( #1 \right)}
\newcommand{\HH}{\mathcal{H}}
\newcommand{\GG}{\mathcal{G}}

\graphicspath{{figures/}} 

\begin{document}

	\title{Flexible rational approximation \\ and its application for matrix functions}
	\date{}
	\author{Nir Sharon\thanks{Corresponding author: nsharon@tauex.tau.ac.il}, Vinesha Peiris, Nadia Sukhorukova,  Julien Ugon}
	\maketitle
	
	\begin{abstract}
		This paper proposes a unique optimization approach for estimating the minimax rational approximation and its application for evaluating matrix functions. Our method enables the extension to generalized rational approximations and has the flexibility of adding constraints. In particular, the latter allows us to control specific properties preferred in matrix function evaluation. For example, in the case of a normal matrix, we can guarantee a bound over the condition number of the matrix, which one needs to invert for evaluating the rational matrix function. We demonstrate the efficiency of our approach for several applications of matrix functions based on direct spectrum filtering.
	\end{abstract}
	



\section{Introduction} 

The term matrix function refers to lifting a scalar function to a matrix domain with matrix values. Such lifting is straightforward for a polynomial since addition and powers are well-defined for square matrices. For nonsmooth functions, the use of polynomials to evaluate matrix functions is limited due to the nature of polynomial approximation. Therefore, rational approximation, which is also well-defined for matrices, introduces an alternative with much-preferred approximation capabilities~\cite{trefethen2019approximation} and comes with the cost of (at least) one inversion. Furthermore, several robust methods for rational approximation have been established in recent years, e.g.,~\cite{filip2018rational, gonnet2013robust, nakatsukasa2018aaa}. Nevertheless, applying rational approximation for evaluating matrix functions remains a challenging task~\cite{nakatsukasa2016computing}.

This paper focuses on an optimization approach for the min-max (uniform) approximation. Perhaps surprisingly, we showed in~\cite{peiris2021generalised} that from the optimization perspective, the min-max problem is more tractable than the corresponding least squares. Moreover, the flexibility of our optimization approach allows us to suggest an improved algorithm, which also includes additional constraints, for example, bounds on the denominator. \ns{In this paper, we demonstrate that this property can be used for improving the optimization algorithms in the case of  evaluating matrix functions.}

 Matrix functions have drawn attention in recent years, see e.g.,~\cite{abdalla2020special}, and proved to be an efficient tool in applications such as reduced order models~\cite{druskin2017multiscale, frommer2008matrix,rankMatFun}, solving ODEs~\cite{li2018exploring}, engineering models~\cite{defez2016solving}, image denoising~\cite{may2016algorithm} and graph neural network~\cite{levie2018cayleynets}, to name a few. 
In this paper, we show that rational approximation is an attractive choice for evaluating matrix functions in the case of normal matrices. In particular, we provide a parameter that controls the \nir{denominator} matrix's conditioning \nir{in our rational approximation} and enables a unique trade-off between the ideal uniform approximation and the well-conditioning of the matrix to be inverted in the rational approximation. We describe the algorithm in detail, prove the theoretical guarantee on the conditioning, and demonstrate numerically many aspects regarding the algorithm behavior. In addition, we present an algorithm application for matrix filtering and validate the advantages of our approach.

We summarize our main contribution as follows:
\begin{enumerate}
	\item 
	We present a flexible optimization framework with an improved algorithm highlighting our ability to incorporate preferred properties for our rational approximation. This algorithm extends the authors' previous work in~\cite{peiris2021generalised}.
	\item
	The algorithm introduces a novel approach to computing approximation with easy add-ons in the form of constraints, providing a framework for tackling general approximation problems.
	\item
	We offer a unique rational approximation for matrix function, explicitly appealing to evaluate nonsmooth or oscillatory functions.
	\item
	We demonstrate the advantage of our matrix function method for several applications involving matrix spectrum slicing.
\end{enumerate}
The paper also supports reproducibility; all source codes and examples are available in an online open repository.\footnote{\url{https://github.com/nirsharon/RationalMatrixFunctions}}

The paper is organized as follows. In Section~\ref{sec:ProblemFormulation}, we set the notation and present the problem and several of its state-of-the-art solutions. In Section~\ref{sec:FlexRat}, we introduce the optimization algorithm and derive some theoretical guarantees regarding our method. In Section~\ref{sec:NumericalExperiments}, we demonstrate the optimization approach numerically and discuss several matrix function applications, including demonstrations via numerical examples. Finally, in Section~\ref{sec:Conclusions}, we discuss the conclusions and future research directions.  

\section{Problem formulation and related background} \label{sec:ProblemFormulation}
We start by formulating the problem we wish to solve. Then, we introduce some additional required notation and background.

\subsection{Evaluating matrix functions}
A matrix function is the result of lifting a scalar function to the domain of a square matrix with square matrix values (of the same order). This paper focuses on the case of real functions of the form $f \colon \mathbb{R} \to \mathbb{R}$ and considers the matrix set consisting of square matrices of size $k \times k$ with real spectrum. 

When $f$ is a polynomial, such lifting is straightforward since addition and powers are well-defined for square matrices.  When $f$ is not a polynomial, there are several standard methods to define the above-mentioned lifting. For example, if $f$ is analytic, having a Taylor expansion whose convergence radius is larger than the spectral radius of a matrix $A$, then the Taylor expansion $f(x) =\sum_{n=0}^\infty \alpha_n x^n$ yields $f(A) = \sum_{n=0}^\infty \alpha_n A^n$. Another elegant definition arises from the Cauchy integral theorem. If $f$ is not analytic, an alternative is to define $f(A)$ on each of the Jordan blocks of $A$ provided that $f \in C^{m-1}$, where $m$ is the size of the largest Jordan block of $A$. For several equivalent definitions and more details, see, e.g., \cite[Chapter 1]{higham2008functions}.

Recently, matrix functions received considerable attention, as many essential challenges in their evaluations have been addressed, e.g., \cite{fasi2018multiprecision, nadukandi2018computing, noferini2017formula, sharon2018evaluating}. The growing number of studies on this topic aims to provide modern tools in various applications. On the classical side, one finds matrix functions in control theory~\cite{kenney1995matrix}. More theoretic topics, where matrix functions are incorporated, \ns{include}  theoretical particle physics and nuclear magnetic resonance, see \cite[Chapter 2]{higham2008functions}. Modern applications include complex network analysis~\cite{benzi2020matrix}, graph convolutional neural networks~\cite{levie2018cayleynets}, \nir{exponential integration for solving large systems of differential equations~\cite{jimenez2020efficient}}, as well as other applications. Moreover, calculating fundamental matrix functions on specially structured matrices, \nir{such as triangular matrices~\cite{parlett1976recurrence}}, opens the door for many new exciting directions; see, e.g., \cite{stotland2018high} and references therein.

Our paper presents a novel algorithm for approximating matrix functions based on rational approximation. The problem we address is defined as follows. Given a real function $f \colon [a,b] \to \mathbb{R}$, and a normal matrix $A$ with all its eigenvalues inside $[a,b]$, construct an approximating matrix $Q$, such that $Q \approx f(A)$. The function is not assumed to be analytic or smooth. In fact, some fascinating cases consist of merely piecewise smooth functions like the sign, square wave, or absolute value functions. It is worth noting that we may consider matrix $A$ with some of its eigenvalues also in the proximity of $[a,b]$ on the complex plane. In particular, if we set, without loss of generality, $a=-1$ and $b=1$, we can allow some extrapolation for eigenvalues in the Bernstein ellipse, that is, the ellipse in the complex plane having focal points at the endpoints, $\pm 1$, see~\cite{webb2012stability}.

\subsection{Uniform rational approximation}

It has been known for several decades~\cite{boyd2004convex, loeb1960} that the optimization problems that appeared in rational and generalized rational approximation in the sense of~\cite{loeb1960} are quasiconvex. 
One of the simplest methods for minimizing quasiconvex functions is the so-called bisection method for quasiconvex optimization~\cite{boyd2004convex}. The primary step in the bisection method is solving convex feasibility problems. Some feasibility problems are challenging, but there are several efficient methods~\cite{Bauschke.ea:2017}. In the case of rational approximation, the feasibility problems we observe in the bisection method can be reduced to solving (large-scale) linear programming problems and, therefore, can be solved efficiently.

Denote by $\ppp{n}$ the space of polynomials of degree at most $n$. The set of type $(m,n)$ rational real functions is defined as $\rrr{m}{n} = \left\lbrace  p/q \mid q\in \ppp{m} , \, \, p \in \ppp{n} \right\rbrace $. When $m=0$, we have that $\rrr{0}{n} = \ppp{n}$. A common parameter choice is, for example, $(m-1,m)$.

Over the years, it was a common belief that the \ns{accuracy} of rational approximation is similar to that of polynomials. In particular, the convergence rates or error bounds are shared between these two families of approximants; see, e.g.,~{\cite[Chapter 23]{trefethen2019approximation}}. While continuous on the interval $[-1,1]$, the absolute value function is not easy to approximate by a polynomial. Specifically, one needs a polynomial of degree $n$ to achieve an error that decays at a rate of $\frac{1}{n} $. This error rate is induced by the smoothness of the function and due to the lack of derivative at the origin. 
\ns{In a seminal paper, Newman~\cite{newman1964rational} showed that the error bound, in the case of rational approximation, is far superior having the order of $\exp(-\sqrt{n})$.}
This result was eventually improved to the asymptotic rate of $8\exp(-\pi\sqrt{n}) $ for the rational minimax~\cite{stahl1993uniform}.

Rational approximations, however, can be problematic. There are various computational challenges here, for example, spurious poles, also known as Froissart doublets. These pole-like points introduce a tiny residue and appear when the degrees are chosen to be too large; see~\cite {beckermann2018rational}. Over the years, various constructions were presented, from the famous Pad{\'e} approximation, which is based on polynomial reproduction, e.g.,~
\cite{brezinski2013pade, gonnet2013robust} to least-squares techniques~\cite{hokanson2018least, van1992parallel}. 
In this paper, we construct rational approximation via the criterion of uniform approximation error.

Uniform rational approximation for a real function $f$ over a closed segment $[a,b]$ is defined as
\begin{equation}\label{eq:uniform_approximation}
\min_{r \in \rrr{m}{n}} \max_{x \in [a,b]} \abs{f(x) - r(x)} .
\end{equation}
The approximant is also known as the minimax approximation. The problem~\eqref{eq:uniform_approximation} is nonconvex; therefore, it is challenging even for the existing tools of modern optimization and nonsmooth calculus. However, there are many interesting observations that can be used to approach this problem.

Several efficient numerical methods are inspired by the celebrated Remez algorithm originally developed for polynomial approximation~\cite{remez1934calcul}. 
Another possibility is to apply one of the tools of modern nonsmooth calculus~\cite{demyanov2014constructive}, which has been successfully applied to free knots polynomial spline approximation problems. 

Two main methods will serve as a benchmark for testing our rational approximation. 
\begin{itemize}
\item 
The first is the adaptive Antoulas--Anderson algorithm~\cite{nakatsukasa2018aaa}, also known as the AAA algorithm. This method combines two approaches. The first approach is the Antoulas Anderson method \cite{antoulas1986scalar} representing rational function in a barycentric manner where the user gives the support points. The second approach is to select the support points in a systematic, greedy fashion to avoid exponential instabilities. This method does not guarantee optimality in any particular norm; however, it bears some advantages, such as the ability to fit many approximation domains and automatically choose the degree of the polynomials while constructing the approximation. Note that an automatic choice of degree in our optimization method can be made, in essence, for a prescribed minimax error. 
\item The second method is the Remez algorithm for rational approximation~\cite{filip2018rational, mayans2006chebyshev}, which solves~\eqref{eq:uniform_approximation} using the equioscillation property. Recent development introduces another usage of this property but from an interpolation point of view~\cite{hofreither2021algorithm}.
\end{itemize}


\subsection{Optimization based on quasiconvexity}

A function $f \colon \mathbb{R}^n \to \mathbb{R}$ is called \textit{quasiconvex} if 
\begin{equation} \label{eqn:quasiconvex_def}
f(\lambda x + (1-\lambda)y) \le \max \{f(x),f(y)\} .
\end{equation}
Equation~\eqref{eqn:quasiconvex_def} is equivalent to saying that the sublevel set of $f$ is convex. Namely, fix $\alpha \in \mathbb{R}$, then
\begin{equation} \label{def:quasiconvex2}
S_\alpha  = \left\lbrace x \in \mathbb{R}^n  \mid   f(x)\le \alpha    \right\rbrace ,
\end{equation}
is a convex set. Quasiconvexity exhibits similar properties to convexity. In particular, the $\sup$ and $\max$ operators preserve quasiconvexity while summation does not.


Quasiconvex problems (with or without linear constraints) can be treated using computational methods developed for quasiconvex optimization. 
In this paper, we use the bisection method for quasiconvex functions. 
The main difficulty of this method is formulating and solving the so-called convex feasibility problems. However, additional linear constraints can be incorporated into the convex feasibility problem in a straightforward manner. 
Since our applications only require additional linear constraints, this method can handle such applications well.
We refer the reader to~\ref{app:bisection}. For further details on this method, see~\cite[Section 4.2.5]{boyd2004convex} and \cite{peiris2021generalised} within the context of rational approximation.

\section{Flexible uniform rational approximation via optimization}\label{sec:FlexRat} 

This section describes in detail the construction of our optimization algorithm. This algorithm is an improved version of the method presented in~\cite{peiris2021generalised}. Here, we focus on a specific instance of rational approximation with further constraints and guarantees.

\subsection{Optimizing a rational approximation}

\subsubsection{Quasiconvexity and the optimization setup}
Our optimization method aims to solve~\eqref{eq:uniform_approximation}. We therefore denote our approximation by
\begin{equation} \label{eqn:rational_form}
    r(x) = \frac{\alpha^{T}\GG(x)}{\beta^{T} \HH(x)}, \quad \text{subject to} \quad  \beta^{T} \HH(x)>0, \quad \HH\in \left(\ppp{m} \right)^{m+1}, \quad \GG\in \left( \ppp{n} \right)^{n+1},
\end{equation}
where $\HH$ and $\GG$ are basis vectorized functions, $\alpha\in \mathbb{R}^{n+1}$ and $\beta\in \mathbb{R}^{m+1}$ are the decision variables, $( \ppp{k} )^{k+1}$, $k\ns{\in}\{n,m\}$  are $(k+1)$-dimensional vectors, whose components are the monomials of degree at most~$k$. The components of $( \ppp{k} )^{k+1},~k\ns{\in}\{n,m\}$ are also called the basis functions. We can use the same formalization to define a generalized rational function: it is enough to replace the basis functions with other Chebyshev systems of the same dimension, for example, exponential functions. In this paper, we focus on polynomials and use as our basis functions the Chebyshev polynomials of the first kind, as defined in~\eqref{eqn:Cheby_polynomial_trig_def}. For more technical details on these polynomials, see~\ref{apn:Cheby_poly}.

At this point, it is worth rewriting~\eqref{eq:uniform_approximation} in terms of~\eqref{eqn:rational_form} to have
\begin{equation} \label{eqn:rational_minimax}
    \min_{\substack{\alpha\in \mathbb{R}^{n+1} \\ \beta\in \mathbb{R}^{m+1}} } \max_{x \in [a,b]} \abs{f(x) - \frac{\alpha^{T}\GG(x)}{\beta^{T} \HH(x)}} ,
\end{equation}
subject to
\begin{equation} \label{eqn:constraint}
    \beta^{T} \HH(x) > 0 , \quad x \in [a,b].
\end{equation}
It has been proved in~\cite[Lemma~2]{Cheney.ea:1964} that the objective function of~\eqref{eqn:rational_minimax} is quasiconvex.\footnote{This important result was not considered as essential for~\cite{Cheney.ea:1964} and therefore was obtained as an intermediate result and went unnoticed for several decades. A simpler proof of this result can be found in~\cite{boyd2004convex, peiris2021generalised}.}. Thus, in light of~\eqref{def:quasiconvex2}, the problem~\eqref{eqn:rational_minimax}-\eqref{eqn:constraint} has an optimal solution for any fixed, large enough (larger than the minimax uniform error), upper bound. 

\subsubsection{Bisection method and convex feasibility problem}
In practice, we solve the problem for a discrete set of points $\{x_i\}_{i=1}^N \subset [a,b]$ where $N \ge m+n+2 $. Namely, we search for a minimal $z$ that solves
\begin{equation} \label{eqn:Opt1}
f(x_i)-\frac{\alpha^{T}\GG(x_i)}{\beta^{T}\HH(x_i)}\le z  \quad  \text{and} \quad 
\frac{\alpha^{T}\GG(x_i)}{\beta^{T} \HH(x_i)}-f(x_i)\le z ,
\end{equation} 
with the constraint
\begin{equation} \label{eqn:Opt3}
\beta^{T} \HH(x_i)>0 .
\end{equation}
To obtain the corresponding sublevel set, fix~$z$  for \eqref{eqn:Opt1}, then the sublevel set is described as \eqref{eqn:Opt1}-\eqref{eqn:Opt3}. Note that when $z$ is fixed, all the constraints are linear, and the convex feasibility problem is reduced to finding a point from this polyhedron. This can be done by solving the following linear programming problem:
$$\min(\theta)$$ 
subject to
\begin{equation} \label{eqn:lp1}
(f(x_i)-z)\beta^{T}H(x_i)-\alpha^{T}\GG(x_i)\le \theta
\end{equation} 
\begin{equation} \label{eqn:lp2}
\alpha^{T}\GG(x_i)-(f(x_i)-z)\beta^{T} \HH(x_i)\le \theta
\end{equation} 
\begin{equation} \label{eqn:lp3}
\beta^{T} \HH(x_i)\geq \delta,
\end{equation}
where $\delta$ is a small positive number.
The feasibility problem is solved if and only if $\theta\leq 0$. 

\begin{remark}
Note that as the number of available points $N$ grows, we obtain more information for approximating the function. However, as $N$ gets larger, so are the required computational resources and runtime of any optimization algorithm that we may employ. Therefore, ideally, we wish to keep $N$ small while extracting enough information to make the rational approximation as accurate as possible. 
\end{remark}

Sampling a function is the subject of numerous studies, e.g.,~\cite{zayed2018advances}, which is beyond the scope of this paper. In practice, one can use a global sampling strategy, for example, equidistant sampling, which is sufficient in many cases. Furthermore, a universal approach for a wide class of common functions is introduced in~\cite{avron2019universal}, which randomly samples the segment of interest according to a non-uniform ideal sampling distribution. Alternatively, when the number of points is required to be small, one may use Chebyshev points (or nodes), which are the roots of the Chebyshev polynomial of the first kind. The Chebyshev polynomials form an orthogonal system~\eqref{eqn:discreteOrtho} (see~\ref{apn:Cheby_poly}), and therefore, they can be utilized as the basis functions ($\HH$ and $\GG$ of~\eqref{eqn:rational_form}). In this case, employing Chebyshev points as sampling points is particularly beneficial to avoid numerical instabilities.


Even though more efficient methods may be applied, we chose the bisection technique for its simplicity and granted convergence to any required tolerance. Several more efficient methods may be applied here, for example, the one from~\cite {pachon2009barycentric}. However, when dealing with matrix functions, the main cost component comes from matrix multiplications and inversion, not accessory computations. Moreover, it is easy to add linear constraints to address a number of specific approximation properties.  
\ns{ For example, if in (\ref{eqn:Opt1}) for a given point $x_i$ one of the constraints (or both) should not exceed a specific positive constant $\varepsilon_i$, that is 
\begin{equation} \label{eqn:Opt11}
f(x_i)-\frac{\alpha^{T}\GG(x_i)}{\beta^{T}\HH(x_i)}\le \varepsilon_i  \quad  \text{and} \quad 
\frac{\alpha^{T}\GG(x_i)}{\beta^{T} \HH(x_i)}-f(x_i)\le \varepsilon_i .
\end{equation}
This additional requirement may be necessarily when some of the grid points require accurate approximations. For example, if one needs a better accuracy at more recent points of a time series or towards the end points of the approximation interval in extrapolation problems. 
}

More work may be required if the additional constraints are convex, quasiaffine (quasilinear), or quasiconvex, but the sets appearing in feasibility problems in the bisection method remain convex since they represent sublevel sets of quasiconvex functions.


\subsection{Denominator Bounds} \label{subsec:denom_bound}

A critical stage in evaluating matrix function via rational approximation of the form $r(x) = p(x)/q(x)$ is to invert the resulting denominator, that is, to calculate $q(A)^{-1}$. Thus, we suggest adding a constraint for bounding~\eqref{eqn:Opt3} as,
\begin{equation} \label{eqn:bounded}
0< \ell \leq \beta^{T} \HH(x_i)\leq u .
\end{equation} 
In the following, we denote the condition number of a matrix $X$ with respect to the matrix $2$-norm by $\cond{X} = \norm{X}_2 \norm{X^{-1}}_2$. In addition, we denote by $\lambda_{\min}(X)$ and $\lambda_{\max}(X)$ the eigenvalues with, respectively, the smallest and largest absolute value of the real matrix $X$. The next lemma shows that bounding the denominator of our approximant from both sides as in~\eqref{eqn:bounded} keeps the condition number of the resulting denominator bounded as well. Note that this lemma can be considered as a direct application of the general result of~\cite{crouzeix2017numerical}.
\begin{lemma} \label{lem:condition_number}
	Let $r=\dfrac{p}{q} \in \rrr{m}{n}$ and assume that
	\begin{equation*} 
	0 < \ell \le \abs{q(x)} \le u , \quad   x \in \Omega . 
	\end{equation*}
Furthermore, assume $A$ is a real normal matrix with all eigenvalues in $\Omega$. Then, 
	\[  \cond{q(A)} \leq u/{\ell} , \] 
\end{lemma}
\begin{proof}
	Denote by $\{\lambda_j\}_{j=1}^k$ the eigenvalues of $A$. Then, the eigenvalues of $q(A)$ are $\{q(\lambda_j)\}_{j=1}^k$ and so,  
	\[ \norm{q(A)}_2 = \max_{\lambda_j} \abs{q(\lambda_j)}  \quad \text{and} \quad  \norm{q(A)^{-1}}_2 = 1/\min_{\lambda_j} \abs{q(\lambda_j)} . \]
Thus, by the bounds on $q$ over $\Omega$ we deduce that
\[  \cond{q(A)} \le  \norm{q(A)}_2  \norm{q(A)^{-1}}_2 \le  u/\ell  . \]
\end{proof}

Lemma~\ref{lem:condition_number} implies that adding bounds on the denominator, when optimizing for the best rational approximation, results in a matrix $q(A)$ with a bounded condition number, regardless of the distribution of eigenvalues of $A$. 


\section{Approximating matrix functions} \label{sec:NumericalExperiments}

The notable advantage of rational approximations over polynomials in scalars motivates us to apply it to matrix functions. Calculating rational matrix functions includes several challenges; for example,~\cite {nakatsukasa2016computing}. However, since matrix polynomials are well understood, the main bottleneck is evaluating the denominator polynomial, which means inverting a matrix polynomial. Two main issues lie in this denominator evaluation. The first is the expensive computational cost. Second, and perhaps more severe in certain situations, is the condition number of the matrix to be inverted. A high condition number indicates the inversion is ill-posed, and the matrix function approximation will contain significant errors. Our general approach for evaluating matrix functions via scalar rational approximation is illustrated in a diagram on Figure~\ref{fig:diagram}.

\begin{minipage}{\textwidth}
\begin{center}
\begin{tikzcd}[row sep=huge,column sep=width("bbbbbbbbbbbbbbbbb")] 
f(x)  \arrow{r}{\text{Scalar approximation}} \arrow[swap]{d}{\text{Matrix function}} & r(x) 
\in \rrr{m}{n}  \arrow{d}{
	\text{Matrix polynomials and one inverse}
} \\
f(A)   & r(A) \arrow{l}{\text{Matrix approximation} }
\end{tikzcd}
 \captionsetup{hypcap=false} 
\captionof{figure}{Approximating matrix function via scalar rational approximation. } \label{fig:diagram}
\end{center}
\medskip
\end{minipage}

In the following, we introduce two applications of our rational approximation for two different tasks of matrix evaluation. The first one is filtering the spectrum of a given matrix, also known as spectrum slicing. The second one is projecting a symmetric matrix to the cone of positive semidefinite matrices. For both applications, the goal is to introduce an alternative solution that does not explicitly recover the spectral structure of the given matrix. Instead, we design and lift a suitable function to obtain the desired matrix function.

The entire source code, including all test parameters, is available at the same GitHub depository as in the previous section. 
We ran the tests over a 2.9 GHz Dual-Core Intel Core i5 processor with 16 GB  RAM using a macOS operation system and Matlab 2017b.

\subsection{Spectrum filtering of a matrix} \label{sssec:filtering_matrix}

The power method is a classical algorithm to extract the leading eigenvalue of a matrix. Its variants and modern versions appear in many applications and areas of science and engineering. When the eigenvalue we wish to recover is not the largest in magnitude, we can use the inverse power method technique and search for the largest eigenvalue of $(A-\alpha I)^{-1}$ where $\alpha$ is a value close to the required eigenvalue, see, e.g.,~\cite{saad2011numerical}. Nevertheless, when the spectrum is unknown, and we wish to preserve or manipulate only a subset of the eigenvalues, we need to search for it extensively. Such direct methods may be costly in terms of computational efforts and generally do not scale well. An alternative approach is to filter the spectrum implicitly by lifting it using a matrix function. This idea was used for spectrum slicing, see, e.g.,~\cite{benner2013computing, li2019eigenvalues}, but can be applied for various broad areas of applications such as improving covariance estimation via shrinkage~\cite{ledoit2020analytical} or optimizing risk factors in financial portfolios~\cite{bouchaud2009financial}, to name a few.

\begin{figure}[ht!]
\centering
         \includegraphics[width=.45\textwidth]{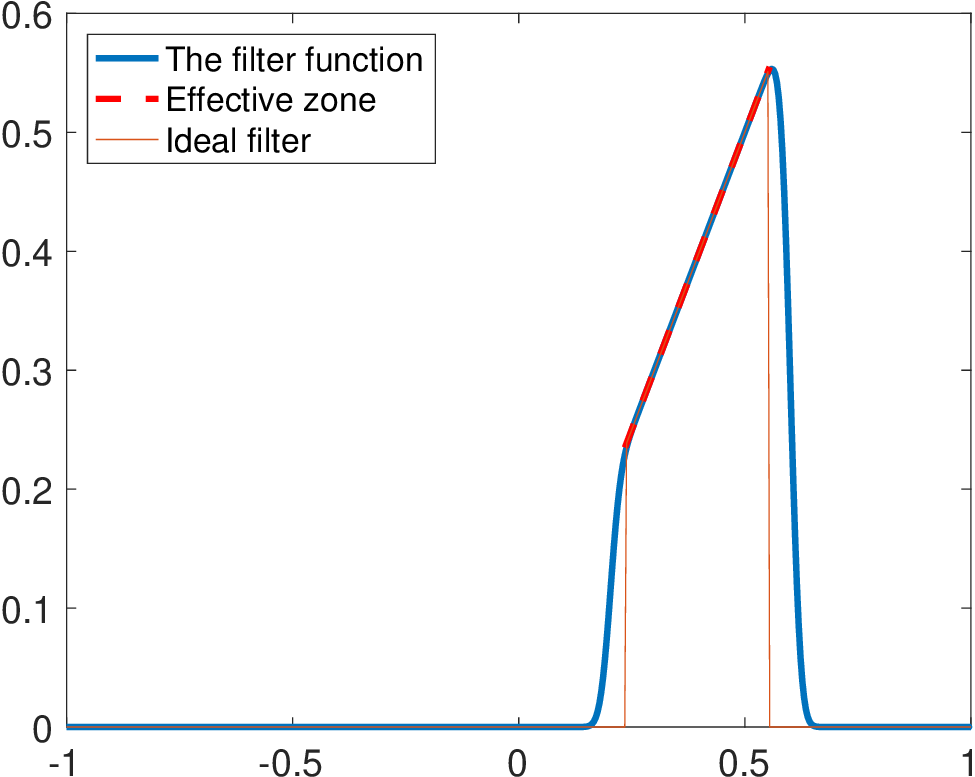}
         \caption{The filter function $F$ of~\eqref{eqn:filter_func}. The effective zone of eigenvalue preservation, $[0.25 , 0.55]$, is marked, together with the corresponding ideal filter function}
         \label{fig:filter_func}
\end{figure}

In the next example, we assume that the eigenvalues of a given matrix are real and bounded in absolute value by $1$. The goal is to retain only the eigenvalues in a certain subsegment of $[-1,1]$. In this example, we choose a segment around $0.4$ at radius $0.15$, that is, $[0.25, 0.55]$. In this case, the ideal filter to be applied on the spectrum is $H(x) = x(h(x-0.25)-h(x-0.55))$ where $h(x)$ is the Heaviside step function. However, this function is not continuous, and so we use a smoother version of it,
\begin{equation} \label{eqn:filter_func}
    F(x) = \frac{x}{2} \left( 1 -  \operatorname{erf}(2( \nicefrac{\abs{x-c}}{r}) - R ) \right). 
\end{equation}
Here, $\operatorname{erf}$ is the Gauss error function, $c = 0.4$ is the center of the segment, R = 0.2 is the width (radius), and $r = 0.05$ determines the rise rate from zero. This function is continuous with a single jump discontinuity at the first derivative (due to the absolute value). The function, the particular segment where we want the eigenvalues to be preserved, and the ideal filter described above are presented in Figure~\ref{fig:filter_func}.

\begin{figure}[ht!]
     \centering
     \begin{subfigure}[t]{0.32\textwidth}
         \includegraphics[width=\textwidth]{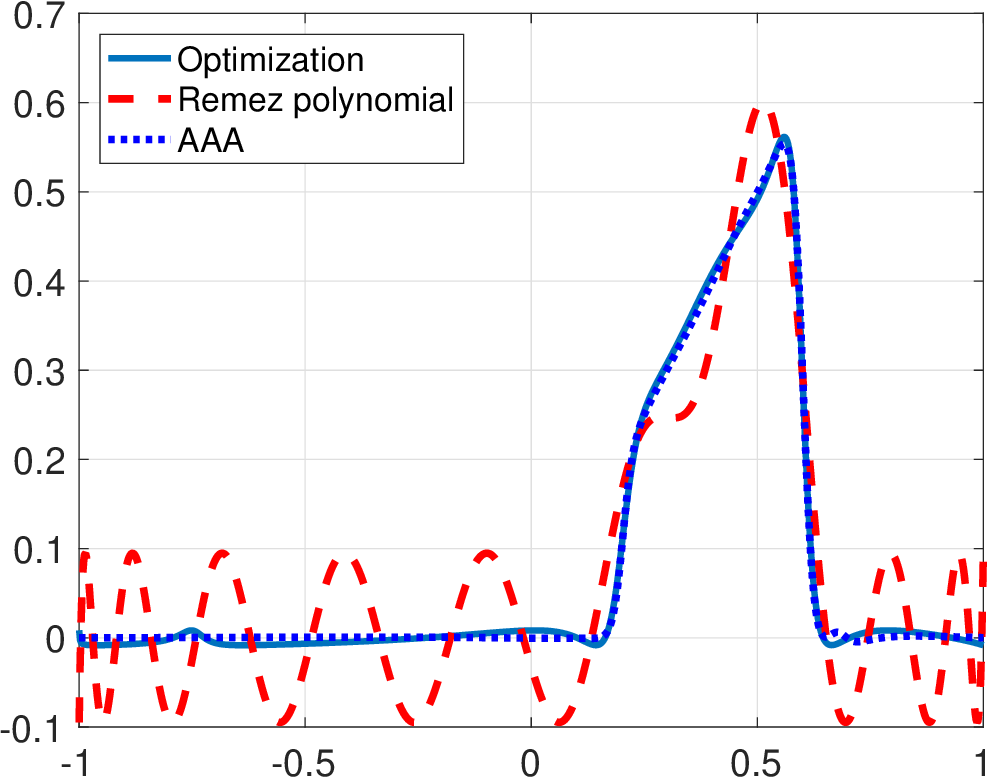}
         \captionsetup{justification=centering}
         \caption{The approximations: $(10,10)$ rational approximations based on our method of optimization and AAA, and a $(20)$ degree minimax polynomial obtained by the Remez algorithm}
         \label{fig:filter_app1}
     \end{subfigure}
     \hfill
     \begin{subfigure}[t]{0.32\textwidth}
         \includegraphics[width=\textwidth]{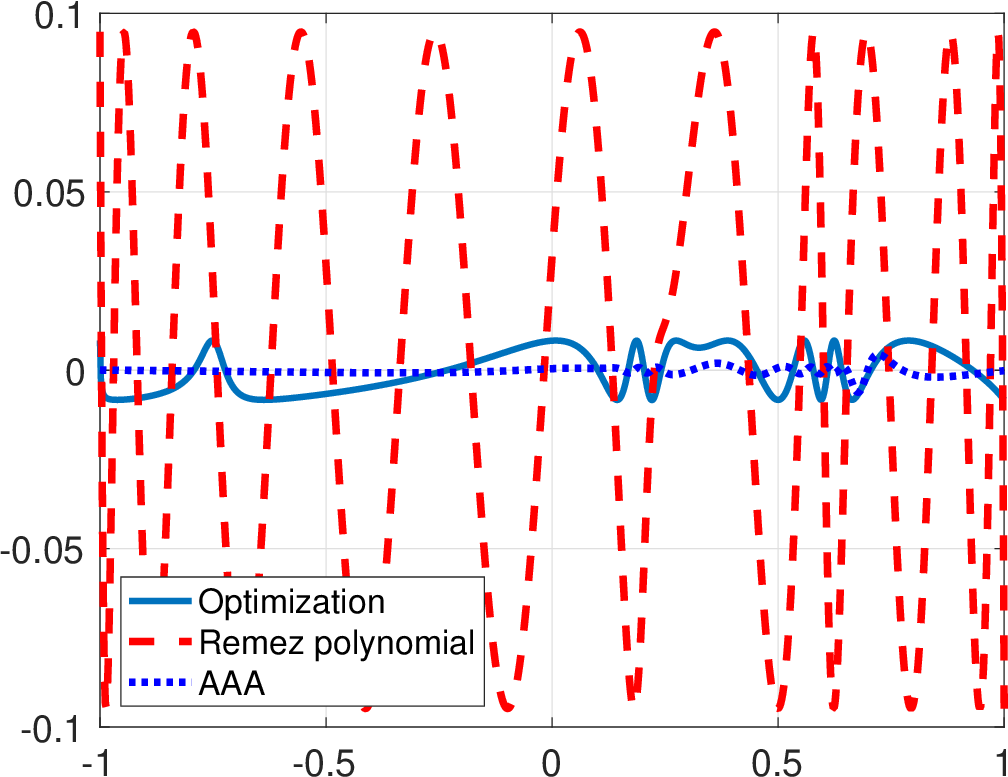}
         \captionsetup{justification=centering}
         \caption{The error rates. In uniform norm, the errors are $0.0083$, $0.0062$, and $0.0948$ for the optimization, AAA, and Remez's polynomial, respectively}
         \label{fig:filter_app2}
     \end{subfigure}
     \hfill
     \begin{subfigure}[t]{0.32\textwidth}
         \includegraphics[width=\textwidth]{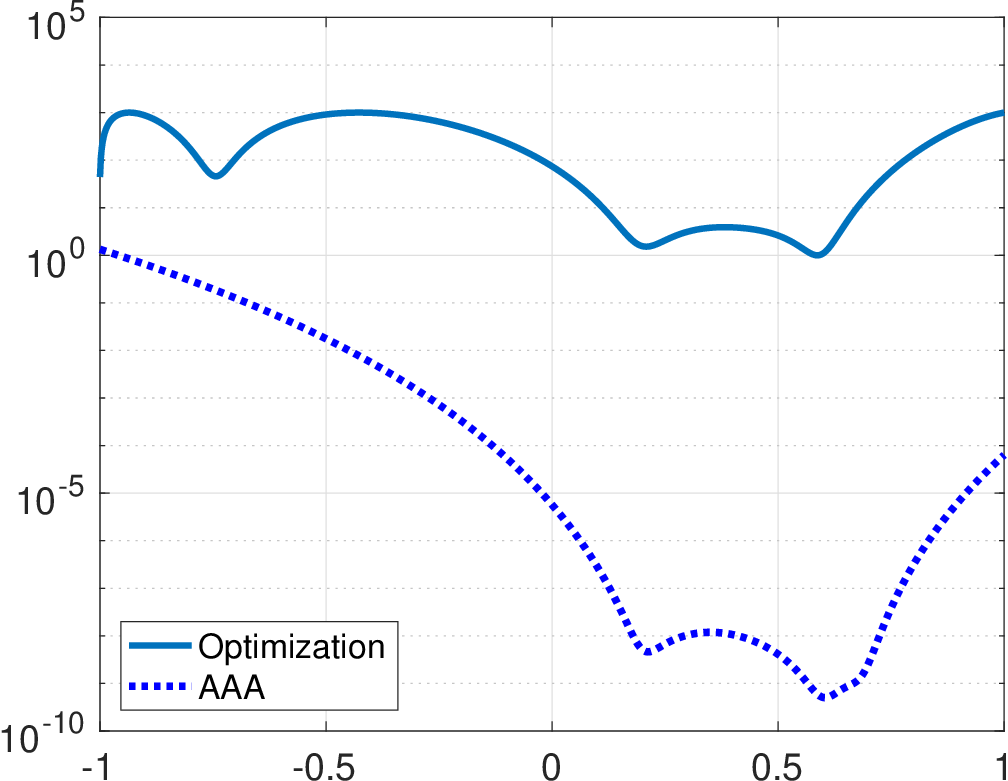}
          \captionsetup{justification=centering}
         \caption{The denominators of the rational approximations. Overall change value $C_r$ of~\eqref{eqn:maxChangeDeno} is $1000$ and $2.6 \times 10^9$ for the optimization and AAA, respectively}
         \label{fig:filter_app3}
     \end{subfigure} 
     \caption{The approximation of the filter function $F$ of~\eqref{eqn:filter_func}}
\label{fig:filter_app}
\end{figure}

We approximate $F$ of~\eqref{eqn:filter_func} using three functions: the output of our optimization, the AAA rational approximation, and the minimax polynomial, as obtained by the Remez algorithm. The two rational approximations are of type $(10,10)$, and the polynomial is of degree $20$. Figure~\ref{fig:filter_app} presents the approximations; in particular, Figure~\ref{fig:filter_app1} shows the functions. In Figure~\ref{fig:filter_app2}, we present the corresponding error rates, which in uniform norm are $0.0083$, $0.0062$, and $0.0948$ for the optimization, AAA, and Remez's polynomial, respectively. For our optimization, we set the upper bound $u$ of~\eqref{eqn:bounded} to be $1000$. The AAA approximation, on the other hand, gives a $C_r$ value of about $2.6 \times 10^9$. The two denominators are given in Figure~\ref{fig:filter_app3}. As we will see next, the change in the denominator of the AAA approximation plays a significant role when applying the rational approximation to matrices.

\begin{figure}[ht!]
	\centering
	\includegraphics[width=.5\textwidth]{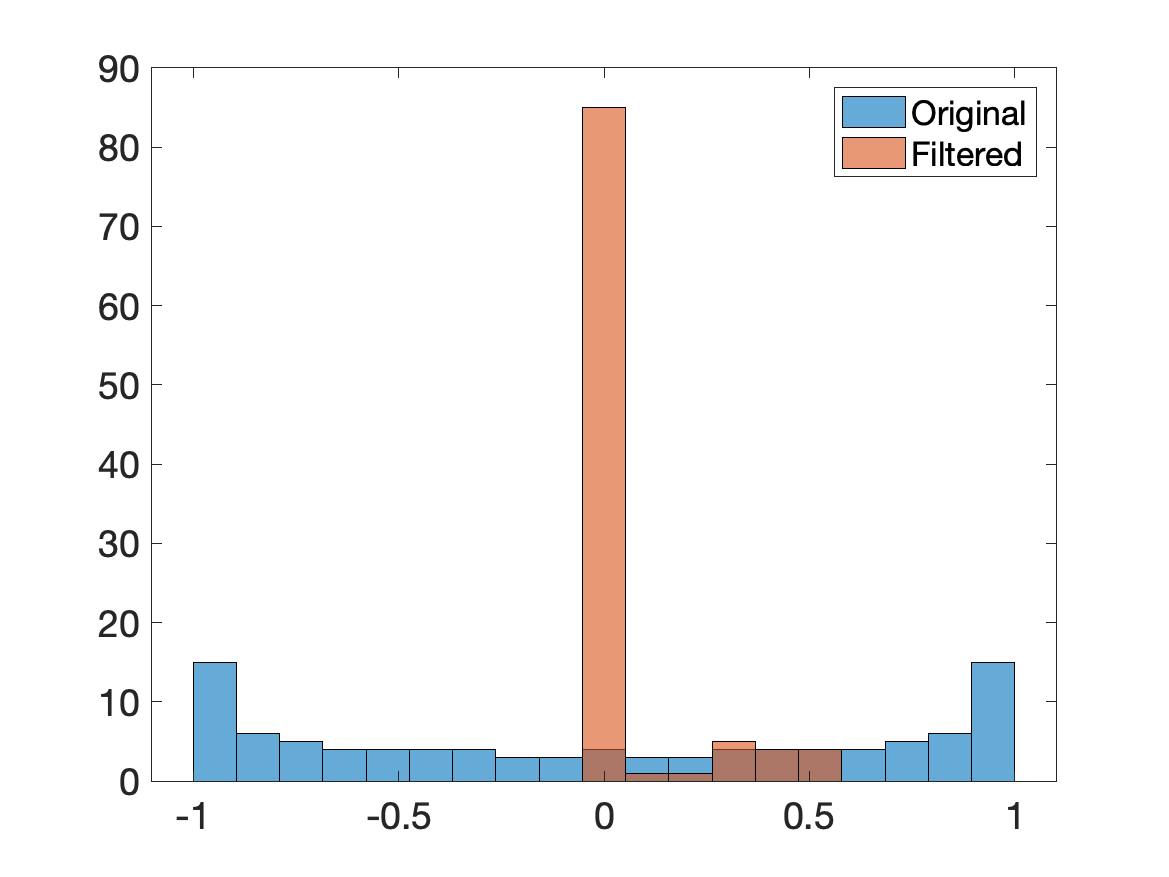}
	\caption{\nir{Algebraic multiplicity of the eigenvalues of the matrix, before and after applying the filter $F$ of~\eqref{eqn:filter_func}.}}
	\label{fig:hist}
\end{figure}

After demonstrating the rational approximations, we apply them to a matrix. In this example, we construct a $100 \times 100$ matrix by randomly selecting an orthogonal matrix $Q$ and using a diagonal matrix $D$ to form $A = QDQ^\ast$. Here, we select $D$ as a diagonal matrix with Chebyshev nodes as its values. Note that $F(A) = QF(D)Q^\ast$ where $F$ is the function given in~\eqref{eqn:filter_func} and $F(D)$ is understood as applying the function on the diagonal element-wise. Each approximation of $F$ will approximate $F(A)$ when applied to $A$. The specific evaluations of the rational polynomials and Remez's polynomial are done using the Horner method, followed by an LU-based inversion for the rational functions. We note that this evaluation method is not optimal as the degree increases, such as, e.g., the Paterson-Stockmeyer~\cite{fasi2019optimality}, but sufficient for our demonstration. \nir{The algebraic multiplicity of the eigenvalues of both $A$ and $F(A)$ are presented in Figure~\ref{fig:hist}}.  

For each approximation $s$ of $F$, we measure its relative error by the Frobenius norm, that is 
\[ \norm{F(A)-s(A)}_{\mathcal{F}}/\norm{F(A)}_{\mathcal{F}} .  \]
The error values were $0.45$ for the minimax polynomial, $0.039$ for our method, and $0.007$ for the more accurate AAA, which was evaluated on the matrix using Horner's rule (evaluating in barycentric form is extremely costly and risky in terms of conditioning). These results were obtained using the double-precision format. So, the above error rates do not reflect the potential nine digits lost by the ill-conditioned denominator of the AAA. Nevertheless, we obtain different results when we repeat the calculations in single precision. While the minimax polynomial and our rational approximation retain their error rates of $0.45$ and $0.039$, the accuracy of AAA is lost. In particular, the relative error reaches a value of $8$, way beyond $100 \%$ of relative error. This means, for example, that the zero function is a better approximation. Note that this scenario is not hypothetical, as most modern GPUs that accelerate computations only support single precision.


\subsection{Applying filtered matrix to a vector} \label{sssec:runtime_vec_mat}

Filtering the spectrum is also applicable in the case of treating the matrix as an operator, that is when we wish to evaluate $f(A)v$ for a certain filter $f$ and vector $v$. Here, we consider a scenario where a rational approximation $r(A)v \approx f(A)v$ is sufficiently accurate, but the evaluation time may be critical. Therefore, for $r(x) = p(x)/q(x)$, we first evaluate the numerator polynomial by adapting to matrix-vector action the Clenshaw’s algorithm, see, e.g.,~\cite[Chapter 5.8]{press2007numerical}. The result is a fast algorithm for evaluating $p(A)v$. Next, we calculate $q(A)$ and solve the linear system $q(A)w = p(A)v $ to obtain $w = q(A)^{-1}p(A)v \approx r(A)v$. In practice, we use Matlab's built-in solver for the latter system, which is LU-based. A faster alternative may be to employ Krylov space methods such as in~\cite{cipolla2021regularization, guttel2013rational}.

\begin{figure}[ht!]
     \centering
     \begin{subfigure}[b]{0.4\textwidth}
         \includegraphics[width=\textwidth]{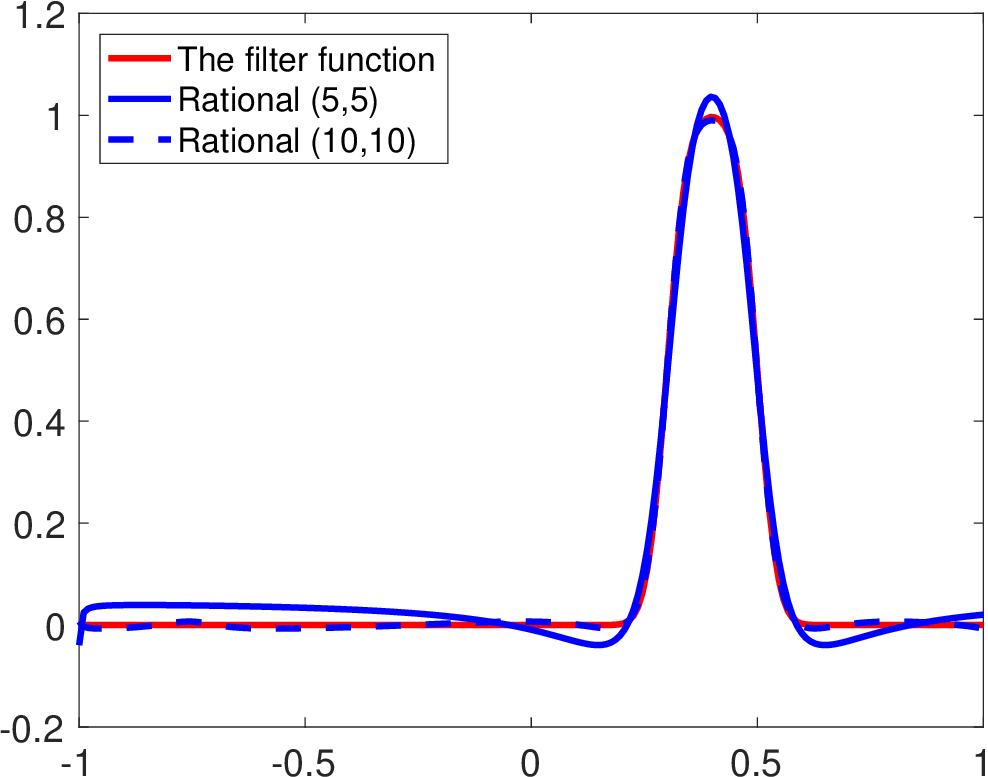}
         \captionsetup{justification=centering}
         \caption{The filter function and its $(5,5)$ and $(10,10)$ rational approximations}
         \label{fig:filter_mat_vec1}
     \end{subfigure}
     \quad
     \begin{subfigure}[b]{0.55\textwidth}
         \includegraphics[width=\textwidth]{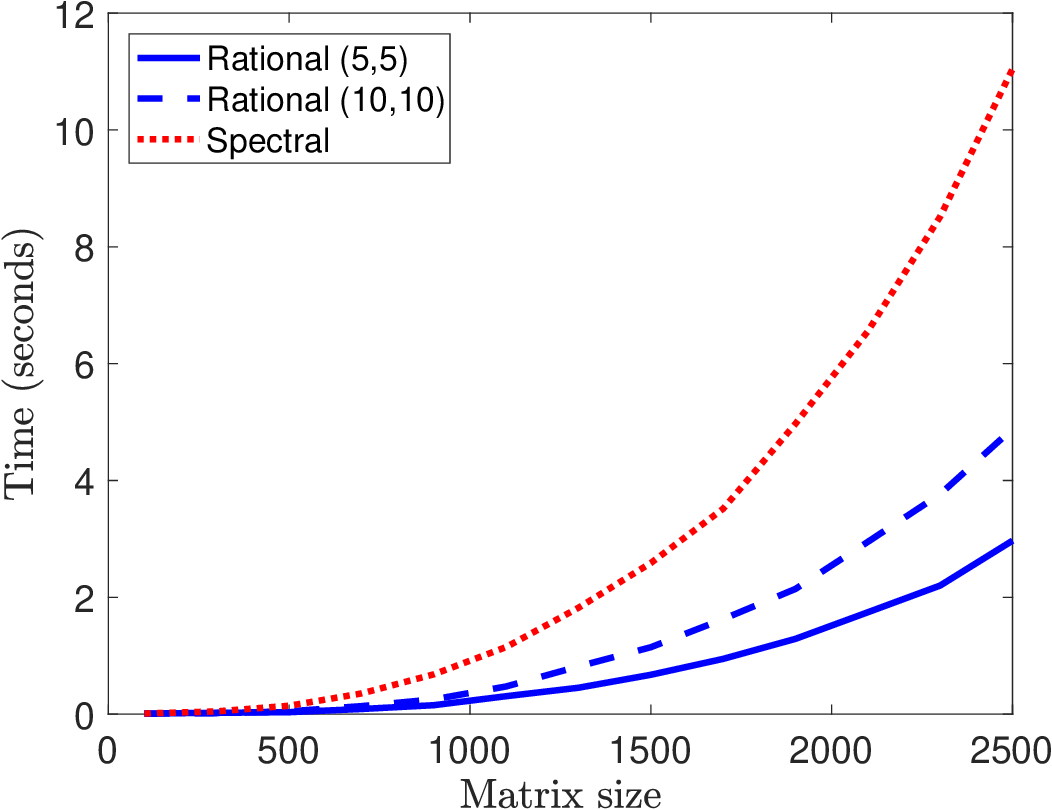}
         \captionsetup{justification=centering}
         \caption{The runtime}
         \label{fig:filter_mat_vec2}
     \end{subfigure}
     \caption{A runtime comparison of evaluating $f(A)v$. The results imply that evaluating the rational approximation is faster than extracting the spectral decomposition, and the difference becomes more significant as the matrix grows}
\label{fig:filter_mat_vec}
\end{figure}

In the following example, we slightly modify the filter F of~\eqref{eqn:filter_func} by setting $R = 0.1$ and $r = 0.1$ and discarding the identity factor. The resulting filter is a sharper bell-shaped function that has small values (roughly) outside the segment $[0.25 , 0.55]$. Therefore, in terms of applying it as a filter for a matrix function, any eigenvalue outside $[0.25 , 0.55]$ will become extremely small. This filter, together with two of its approximations of type $(5,5)$ and type $(10,10)$, are depicted in Figure~\ref{fig:filter_mat_vec1}. The rational functions are calculated with an upper bound of $u=1000$, and their approximation errors, in the uniform norm, are $0.0395$ and $0.0069$ for the $(5,5)$ and $(10,10)$ approximations, respectively. The matrices are of the form $A = QDQ^\ast\in \mathbb{R}^{k \times k}$ where the size $k$ varies from $100$ to $2500$, $Q$ is a randomly selected orthogonal matrix, and $D$ is a diagonal matrix with eigenvalues drawn uniformly at random over $[-1,1]$. As a benchmark, we calculate the spectral decomposition of $A$, using Matlab's \textit{eig} function, and then evaluate the filtered vector as $Q(F(D)(Q^\ast v))$, that is, by three matrix-vector operations. The accuracy is measured in relative error $ \norm{F(A)v-r(A)v}/\norm{F(A)v} $, and while the spectral decomposition presents errors of about $10^{-10}$, the rational functions obtain relative error of about $1\%$ and $5\%$, for the $(5,5)$ and $(10,10)$ approximations, respectively. Nevertheless, the benefit appears in the form of a better runtime, as presented in Figure~\ref{fig:filter_mat_vec2}. The runtime is measured in seconds and averaged over $10$ repeated trials. The outcome suggests that our method is faster, and the gap increases as the size of the matrix grows.

\subsection{From symmetric matrices to positive semidefinite matrices} \label{sssec:projection}

While symmetric matrices have real eigenvalues, these eigenvalues are not necessarily nonnegative.
The problem of projecting a symmetric matrix to the cone of positive semidefinite matrices, that is, finding the nearest positive semidefinite matrix to a symmetric matrix in the spectral norm appears in many studies, see, e.g.,~\cite{boyd2005least}, and applications, for example, in finance industry~\cite{higham2002computing} and risk management~\cite{cutajar2017actuarial}. In its plain form, this problem has a straightforward solution: calculating the spectral decomposition and then zeroing out the negative eigenvalues. This procedure solves projection with respect to the squared Frobenius norm; see, e.g.,~\cite[Chapter 8]{boyd2004convex}. Nevertheless, significant time is required to do so when the matrix size increases. Thus, different methods were proposed over time, e.g.,~\cite{higham2002accuracy, kressner2017fast}, including a spectral projectors analysis~\cite{benzi2013decay}. Here, we provide another alternative approach; we propose to apply the function 
\begin{equation} \label{eqn:relu_fun}
    f(x) = \max\{0,x\}, 
\end{equation}
also known as the ReLU (Rectified Linear Unit) function, to the symmetric matrix. The following examples illustrate this matrix function technique and compare it with the textbook solution based on the native Matlab implementation of spectral decomposition. We demonstrate the flexibility of our method, which now includes positivity as a constraint. We then show that the matrix function approach yields an efficient algorithm to solve the above problem, assuming that a moderate accuracy level is sufficient.

\begin{figure}[ht]
     \centering
     \begin{subfigure}[t]{0.32\textwidth}
         \includegraphics[width=\textwidth]{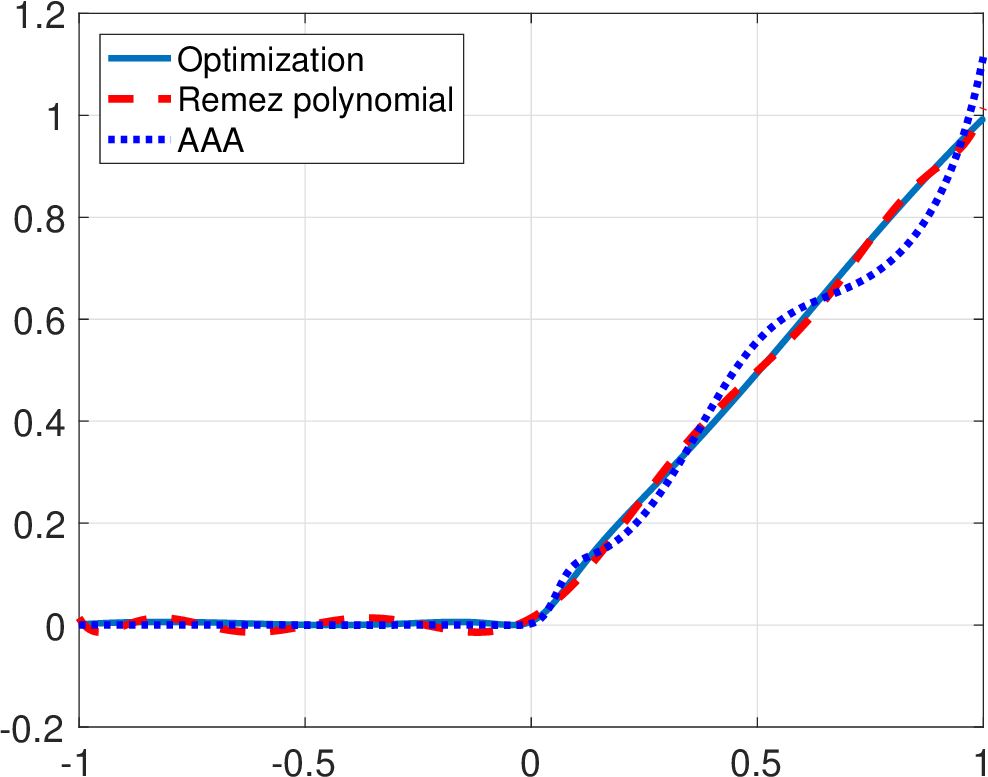}
         \captionsetup{justification=centering}
         \caption{The approximations to the ReLU function: two $(5,5)$ rational approximations based on our method of optimization and of AAA, and a $(10)$ degree minimax polynomial obtained by the Remez algorithm}
         \label{fig:spd_proj1}
     \end{subfigure}
     \hfill
     \begin{subfigure}[t]{0.32\textwidth}
         \includegraphics[width=\textwidth]{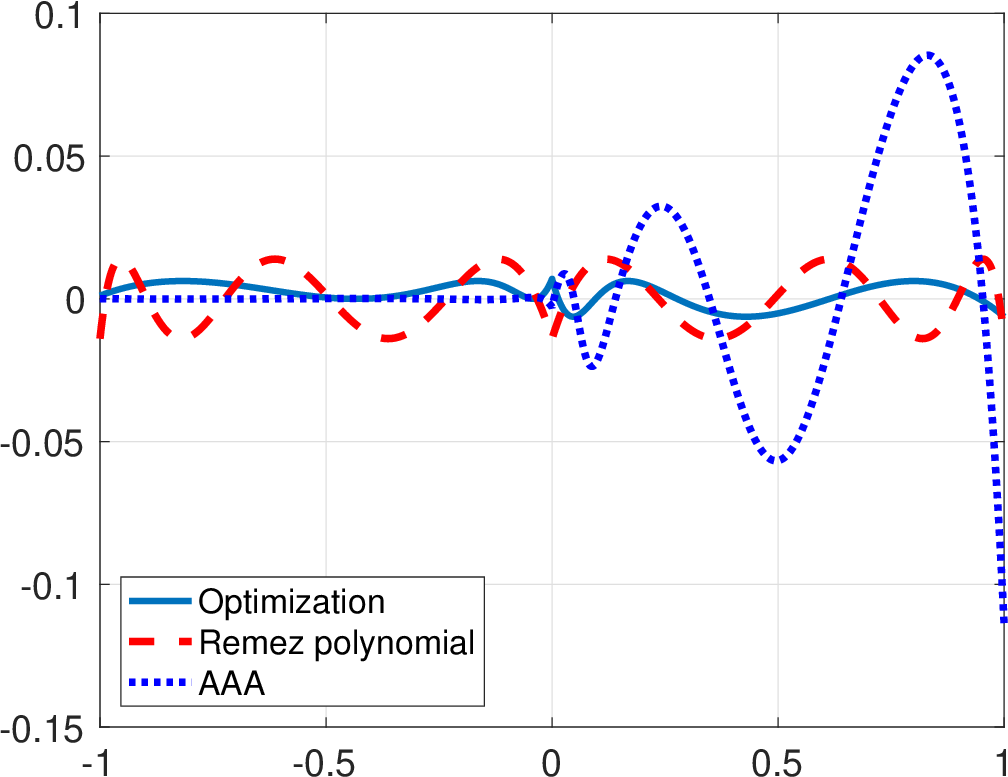}
         \captionsetup{justification=centering}
         \caption{The error rates. In uniform norm, the errors are $0.007$, $0.013$, and $0.114$ for the optimization, AAA, and Remez's polynomial, respectively. Note that the AAA is much accurate on the segment $[-1,0]$ where it presents a uniform error of $0.0017$}
         \label{fig:spd_proj2}
     \end{subfigure}
     \hfill
     \begin{subfigure}[t]{0.32\textwidth}
         \includegraphics[width=\textwidth]{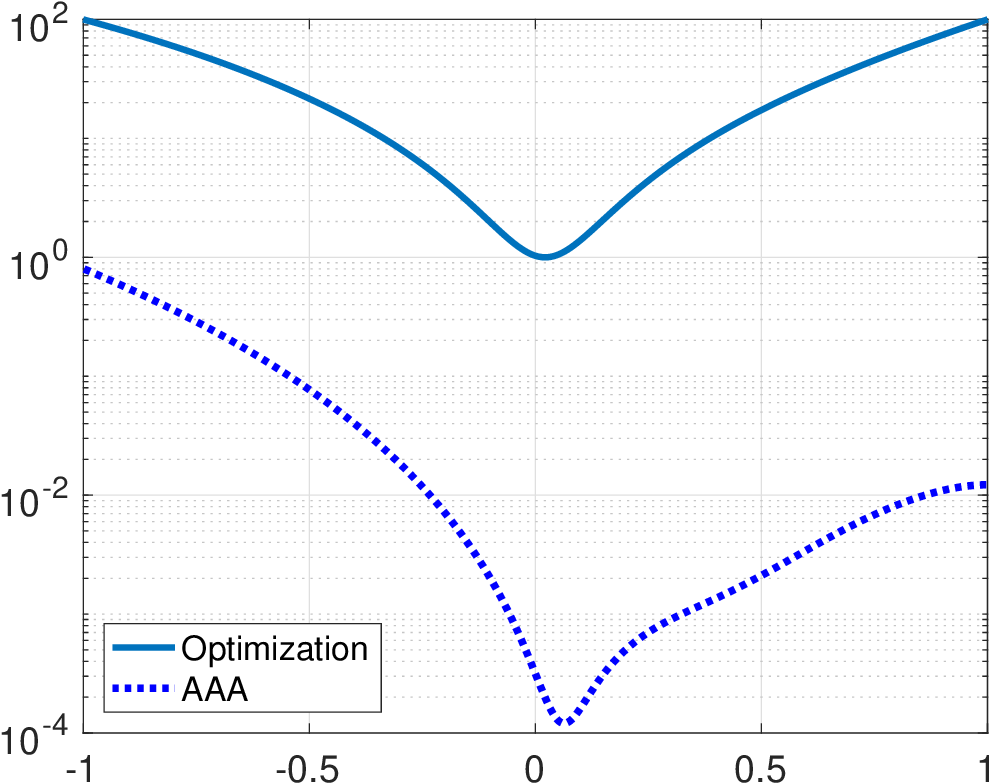}
          \captionsetup{justification=centering}
         \caption{The denominators of the rational approximations. Overall change value $C_r$ of~\eqref{eqn:maxChangeDeno} is $100$ and $6603$ for the optimization and AAA, respectively}
         \label{fig:spd_proj3}
     \end{subfigure} 
     \caption{Approximations of the ReLU function $\max\{0,x\}$}
\label{fig:spd_proj}
\end{figure}

As in the previous examples, we start with approximating the function, which is, in this case, $\max\{0,x\}$. One caveat of polynomial and rational approximations is that the approximant oscillates, so the required non-negativity is hard to achieve. Solutions like using Fej\'er kernel or damping the polynomial coefficients with suitable weights usually result in a slow convergence and lower quality approximation, see, e.g.,~\cite{saad2011numerical}. In our method, we add a positivity constraint to the numerator, merely posing a linear constraint to the optimization problem we solve. The additional restriction reduces the search space, but the resulting accuracy turns out to be comparable in magnitude. In particular, the positivity constraint slightly increases the uniform error from $0.0055$ to $0.007$ for a (5,5) type rational function with a denominator upper bound of $u=100$. The rational approximation, together with the best uniform polynomial and the AAA rational approximation, appear in Figure~\ref{fig:spd_proj}. The figure shows how our rational approximation remains positive and retains the lowest uniform error. In addition, we show how the denominator of the AAA approximation varies across almost four orders of magnitude while our denominator keeps bounded between $1$ and $100$.

We also measure the maximum absolute change in the denominator values of a rational function $r(x) = p(x)/q(x)$ over $[a,b]$, denoted by
 \begin{equation} \label{eqn:maxChangeDeno}
    C_r = \max_{x \in [a,b]} \abs{q(x)} / \min_{x \in [a,b]} \abs{q(x)} .
 \end{equation}
 This value indicates how high the condition number of $q(A)$ may be for a given matrix $A$ with eigenvalues inside $[a,b]$.

\begin{figure}[ht]
\begin{minipage}[b]{0.48\linewidth}
\centering
\includegraphics[width=\textwidth]{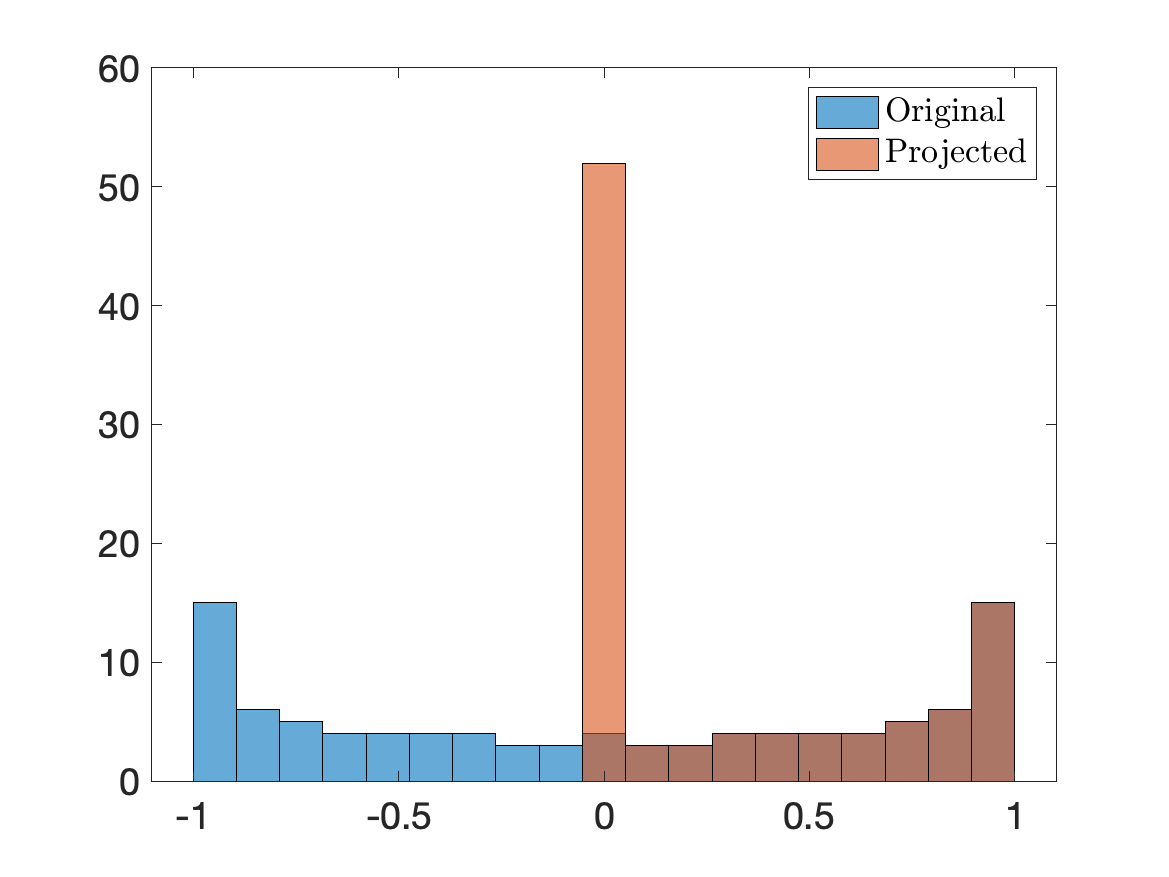}
\caption{\vin{Algebraic multiplicity of the eigenvalues of the symmetric matrix, before and after projection}}
\label{fig:spd_proj_hist}
\end{minipage}
\quad
\begin{minipage}[b]{0.48\linewidth}
\centering
\includegraphics[width=\textwidth]{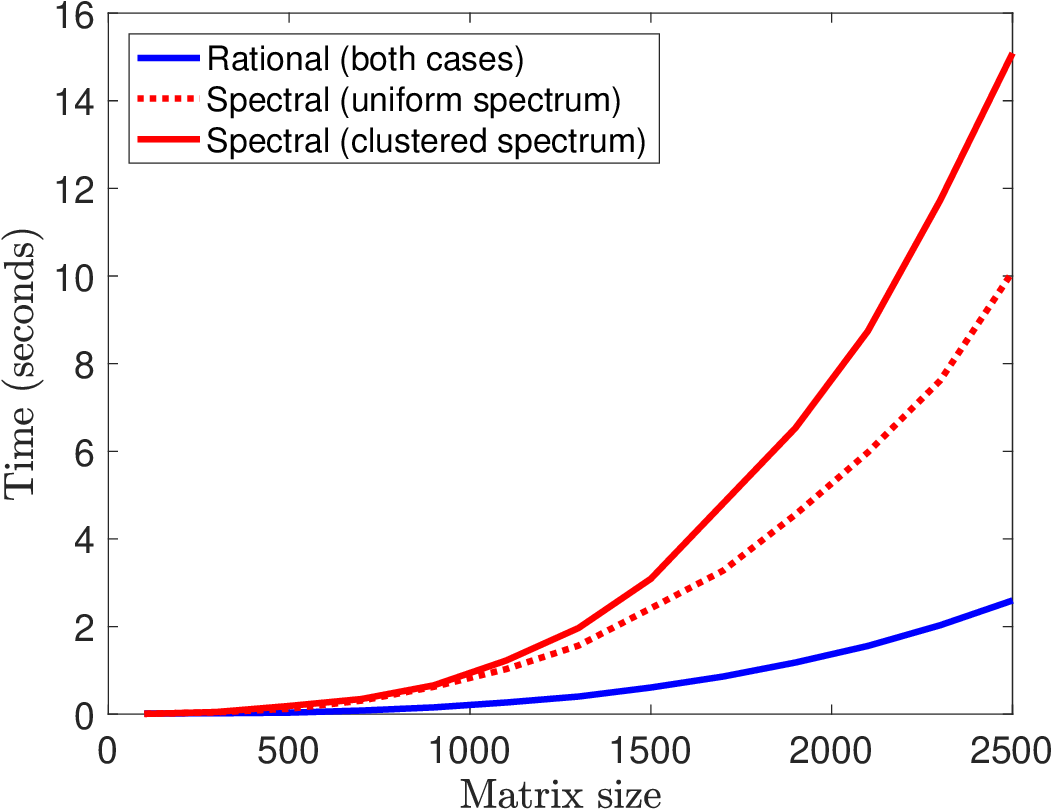}
         \caption{The time comparison of projection onto the cone of SPD matrices}
         \label{fig:projection}
\end{minipage}
\end{figure}

We propose to compute the projection of a symmetric matrix $A$ to the cone of positive semidefinite matrices as the matrix function $r(A)$ where $r(x)\approx \max\{0,x\}$. Similar to Section~\ref{sssec:filtering_matrix}, our first test uses a symmetric matrix of size $100 \times 100$ with its eigenvalues spread uniformly equidistantly over the segment. We apply the $(5,5)$ non-negative rational approximation, as presented in Figure~\ref{fig:spd_proj}, over the matrix. Then, we plot a histogram \nir{representing the algebraic multiplicity of the eigenvalues} of the matrix and its projected version. This histogram is given in Figure~\ref{fig:spd_proj_hist} and shows the effect of projection as all eigenvalues are non-negative.

The second test is a time comparison, similar to Section~\ref{sssec:runtime_vec_mat}. In particular, we compare the evaluating of $f(A)v$ where $f$ is the ReLU function of~\eqref{eqn:relu_fun} and $A$ and $v$ are a given matrix and a vector. The procedure is the same as in Section~\ref{sssec:runtime_vec_mat}, so we omit it for brevity. Nevertheless, in this example, we run the test over two types of matrices. The first has eigenvalues that were uniformly drawn from $[-1,1]$. The second type of matrices has eigenvalues grouped into two narrow segments around $\pm 0.3$. Such clustered spectra are the Achilles' heel of many standard spectral decomposition algorithms. Indeed, as seen in Figure~\ref{fig:projection}, the runtime as a function of the matrix size increases rapidly for these type of matrices compared to the first type. On the other hand, our matrix function method gives similar runtime performances in both cases. 

\section{Conclusions and future research directions}\label{sec:Conclusions}

This paper demonstrates that matrix function lifting, frequently appearing in practical applications, can be tackled through rational and generalized rational approximation of the original function. This approach combines a high level of approximation accuracy and the simplicity of the computational procedures. Another essential advantage of our method is that we may naturally add more constraints to the model without decreasing the method's efficiency. In particular, this advantage is valid when the added constraints are linear. Furthermore, this extension gives rise to some improvements, particularly the ability to control the conditioning number of the matrices. Finally, the numerical experiments demonstrate the efficiency of our method. 

Our future research directions include the possibility of having non-linear constraints as well. The first step is to extend it to the case of quasiaffine (quasilinear) additional constraints since the corresponding feasible sets remain polyhedra.

\section*{Acknowledgment} Nir Sharon is partially supported by the NSF-BSF award 2019752 NS and the DFG award 514588180. Vinesha Peiris, Nadezda Sukhorukova, and Julien Ugon are supported by the Australian Research Council (ARC), Solving hard Chebyshev approximation problems through nonsmooth analysis (Discovery Project DP180100602).

\bibliographystyle{plain}
\bibliography{bibiliography}

\appendix

\section{Bisection method for quasiconvex optimization} \label{app:bisection}

The bisection method is a simple and reliable quasiconvex optimization method, e.g.,~\cite [Section~4.2.5]{boyd2004convex}. In all the examples in this paper, the corresponding convex feasibility problems (discrete case) can be reduced to solving linear programming problems when the space is finite, while for the continuous case, the convex feasibility problems are linear semi-infinite. 

Consider a quasiconvex optimization problem with linear constraints:
\begin{equation}\label{eq:quasiconvexObjFun}
    {\rm minimize}~ \phi(x)
\end{equation}
subject to
\begin{equation}\label{eq:quasiconvexConstraints}
    M x\leq b,
\end{equation}
where $\phi(x)$ is a bounded quasiconvex function, $x\in X\subset\mathbb{R}^d$, $X$ is a convex set. Constraints~(\ref{eq:quasiconvexConstraints}) are linear equations and inequalities. The set, let us call it $C$, of feasible points satisfying these constraints is convex; that is, the set $C$ is nonempty. Fix~$z$, the sublevel set $\phi_z=\{z:\phi(x)\leq z\}$ is convex, since~$\phi$ is quasiconvex.

The set $C_z=C\cap \phi_z$ is convex as the intersection of two convex sets. Assume that problem~(\ref{eq:quasiconvexObjFun})-(\ref{eq:quasiconvexConstraints}) is feasible, that is, constraints~(\ref{eq:quasiconvexConstraints}) are consistent. Then, 
\begin{itemize}
    \item 
if~$C_z$ is empty, then the optimal value for~(\ref{eq:quasiconvexObjFun})-(\ref{eq:quasiconvexConstraints}) is strictly higher than~$z$; 
\item if~$C_z$ is not empty, then the optimal value for~(\ref{eq:quasiconvexObjFun})-(\ref{eq:quasiconvexConstraints}) does not exceed~$z$.    
\end{itemize}
The procedure of checking whether set~$C_z$ is empty or not is called convex feasibility. 

Note that feasibility problems aim at finding a feasible point in a given set of constraints, and therefore this point is not necessarily an optimal solution for some objective function associated with this set of constraints, but in most cases, feasibility problems can be formulated as optimization problems. 

It was noticed a long time ago that rational and generalized rational approximation problems could be solved by fixing the level set of the objective function at some value and then solving linear programming problems, see for example~\cite{meinardus2012approximation, ralston1965rational}. 
Moreover, some earlier developed methods are implementations of bisection method for quasiconvex functions~\cite{PeirisSukhBisecIneq}. It was demonstrated in \cite{rubinov2005abstract} that the sublevel sets of quasiaffine functions are half-spaces. Therefore the constraint sets in the corresponding optimization problems are polyhedra, but there is no approach for constructing these polyhedra for general quasiaffine constraints. However, in the case of rational and generalized rational approximation, we know how to construct the polyhedra, and the idea comes from earlier works. More details are in section~\ref{sec:FlexRat}. This problem is a beautiful example of interconnections between modern optimization techniques and classical approximation results obtained in the 60-70s of the twentieth century.  

\begin{algorithm}
	\caption{Bisection algorithm for quasiconvex optimization with linear constraints} \label{alg:bisection}
	\begin{algorithmic}
			\ENSURE Maximal deviation $z$ (within $\varepsilon$ precision). 
			Start with given precision $\varepsilon>0$.
			\STATE set $l \leftarrow 0$
			\STATE set $u$
			\STATE $z \leftarrow (u+l)/2$
			\WHILE{$u-l \leq \varepsilon$}
			\STATE {Check if $C_z$ is empty.}
			\IF{feasible solution exists ($C_z$ is not empty)}
			\STATE $u \leftarrow z$
			\ELSE
			\STATE $l \leftarrow z$
			\ENDIF
			\STATE update $z \leftarrow (u+l)/2$
			\ENDWHILE 
		\end{algorithmic}
\end{algorithm}

The bisection method which we use is given in Algorithm~\ref{alg:bisection}. It is essential for Algorithm~\ref{alg:bisection} to start with tight upper bound~$u$ and lower bound~$l$. Since our quasiconvex function corresponds to the maximal absolute deviation between original function~$f$ on $[a,b]$ (multivariate generalizations are also possible) and its approximation, the lower bound $l$ can be set as zero (that is, $l \leftarrow 0$). The upper bound~$u$ can be set as $u=\frac{1}{2}(\max f(x)-\min f(x))$.

Note that at each iteration of the bisection algorithm, the length of the search interval is halved, meaning that it takes at most $\log_2(\nicefrac{L}{\varepsilon})$ steps to reach a solution within the desired accuracy $\varepsilon$. When approximating the function over a discrete set of points, the convex feasibility problem can be solved with a polynomial-time algorithm for linear programming. 

\section{Chebyshev polynomials} \label{apn:Cheby_poly}

Chebyshev polynomials of the first kind of degree $n$ are defined as
\begin{equation} \label{eqn:Cheby_polynomial_trig_def}
T_n(x) = \cos \left( n \arccos(x) \right) , \quad x \in [-1,1] , \quad n = 0,1,2,\ldots
\end{equation}
These polynomials are solutions of the Sturm-Liouville ordinary differential equation
\begin{equation} \label{eqn:ODE}
(1-x^2)y^{\prime \prime } - xy^\prime +n^2y=0
\end{equation}
and satisfy the three-term recursion
\begin{equation} \label{eqn:three_term_recursion}
T_n(x) = 2xT_{n-1}(x) - T_{n-2}(x) , \quad n=2,3,\ldots
\end{equation}
with $T_0(x)=1$ and $T_1(x) = x$. Therefore, Chebyshev polynomials form an orthogonal basis for $L_2([-1,1])$ with respect to the inner product
\begin{equation}  \label{eqn:inner_product}
\left\langle  f,g \right\rangle_T =  \frac{2}{\pi} \int_{-1}^1 \frac{f(t)g(t)}{\sqrt{1-t^2}} dt .
\end{equation}
The Chebyshev expansion of a function $f$ with a finite norm with respect to \eqref{eqn:inner_product} is 
\begin{equation} \label{eqn:chebyshev_series_scalars}
f(x) \sim \sumd{\infty} \alpha_n[f] T_n(x) , \quad \alpha_n[f] = \left\langle f, T_n \right\rangle_T ,
\end{equation}
where the dashed sum \resizebox{.33in}{!}{$\left( \sideset{}{'}\sum \right)$} denotes that the first term is halved. The truncated Chebyshev expansion is defined as 
\begin{equation} \label{eqn:truncated_chebyshev}
\mathcal{S}_N(f)(x)  = \sumd{N} \alpha_n[f]T_n(x)  .
\end{equation}
$\mathcal{S}_N(f)(x)$ is a polynomial approximation of $f$ which is the best least squares approximation with respect to the induced norm $\norm{f}_T = \sqrt{\left\langle  f,f \right\rangle_T}$. Remarkably, this least squares approximation is close to the best minimax polynomial approximation, and in particular the following was established by Bernstein \cite{bernstein1918quelques}
\[ \norm{f(x) - S_N(f)(x)  }_{\infty} \le  \Lambda_N \norm{f(x) - p^\ast_N(x)  }_{\infty} , \quad f \in C([-1,1]) ,\]
where $p^\ast_N$ is the unique best minimax polynomial approximation of degree $N$, and Lebesgue constant $\Lambda_N$ behaves asymptotically as $\log(N)$ for large $N$. For the exact value, see~\cite{powell1967maximum} (and in particular, it is less than $6$ for $N<1000$). 

The Chebyshev polynomials satisfy a discrete orthogonality relation as well as the continuous one (with respect to~\eqref{eqn:inner_product}). Specifically, denote by $z_{n}$, $n=1,\ldots, N$ the zeros of $T_{N}$. Then, for any $ 0 \le \ell \neq k<N$ we have
\begin{equation} \label{eqn:discreteOrtho}
    \sum_{n=1}^{N} T_{\ell}(z_{n}) T_{k}(z_{n})=0 .
\end{equation}
The above is essential in calculating $\alpha_n$ of~\eqref{eqn:truncated_chebyshev} and plays a significant role in stabilizing our optimization that uses the representation~\eqref{eqn:rational_form}.
For more information, we refer the interested reader to~\cite[Chapter 5]{press2007numerical}.

\end{document}